\documentclass{article}
\usepackage{amssymb}

\input{tcilatex}
\begin{document}

\begin{center}
\textbf{The Cauchy problem for operator-Boussinesq equations}

\textbf{Veli\ B. Shakhmurov}

Department of Mechanical Engineering, Okan University, Akfirat, Tuzla 34959
Istanbul, Turkey, telephone: +90 216 6771630, fax: +90 216 6771486, E-mail:
veli.sahmurov@okan.edu.tr;

A\textbf{bstract}
\end{center}

In this paper, the existence and uniqueness of solution of the Cauchy
problem for abstract Boussinesq equation is obtained. By applying this
result, the Cauchy problem for systems of Boussinesq equations of finite or
infinite orders are studied.

\textbf{Key Word:}$\mathbb{\ \ }$Boussinesq equations\textbf{, }Semigroups
of operators, Hyperbolic-operator equations; cosine operator functions,
Operator-valued multipliers

\begin{center}
\bigskip\ \ \textbf{AMS: 35Lxx, 35Mxx, 47Lxx, 47Axx }

\textbf{1}. \textbf{Introduction}
\end{center}

The subject of this paper is to study the local existence and uniqueness of
solution of the Cauchy problem for the following Boussinesq-operator equation%
\begin{equation}
u_{tt}-Lu_{tt}+Au=f\left( u\right) ,\text{ }x\in R^{n},\text{ }t\in \left(
0,T\right) ,  \tag{1.1}
\end{equation}%
\begin{equation}
u\left( x,0\right) =\varphi \left( x\right) ,\text{ }u_{t}\left( x,0\right)
=\psi \left( x\right) ,  \tag{1.2}
\end{equation}%
where $A$ is a linear operator in a Banach space $E$, $u(x,t)$ is the $E$%
-valued unknown function, $f(u)$ is the given nonlinear function, $\varphi
\left( x\right) $ and $\psi \left( x\right) $ are the given initial value
functions, subscript $t$ indicates the partial derivative with respect to $t$%
, $n$ is the dimension of space variable $x$ and $L$ is an elliptic operator
in $R^{n}$ with constant coefficients$.$ Since the Banach space $E$ is
arbitrary and $A$ is a possible linear operator, by choosing $E$ and $A$ we
can obtain numerous classes of generalized Boussinesq type equations which
occur in a wide variety of physical systems, such as in the propagation of
longitudinal deformation waves in an elastic rod, hydro-dynamical process in
plasma, in materials science which describe spinodal decomposition and in
the absence of mechanical stresses (see $\left[ 1-4\right] $ ). For example,
if we choose $E=\mathbb{C},$ $L=\Delta $ and $A=-\Delta $ we obtain the
scalar Cauchy problem for generalized Boussinesq type equation%
\begin{equation}
u_{tt}-\Delta u_{tt}-\Delta u=f\left( u\right) ,\text{ }x\in R^{n},\text{ }%
t\in \left( 0,T\right) ,  \tag{1.3}
\end{equation}%
\begin{equation}
u\left( x,0\right) =\varphi \left( x\right) ,\text{ }u_{t}\left( x,0\right)
=\psi \left( x\right) .  \tag{1.4}
\end{equation}%
The equation $(1.3)$ arises in different situations (see $[1,2]$). For
example, for $n=1$ it describes a limit of a one-dimensional nonlinear
lattice $\left[ 3\right] $, shallow-water waves $\left[ 4,5\right] $ and the
propagation of longitudinal deformation waves in an elastic rod $\left[ 6%
\right] $. Rosenau $\left[ 7\right] $ derived the equations governing
dynamics of one, two and three-dimensional lattices. One of those equations
is $\left( 1.3\right) $. In $\left[ 8\right] $, $\left[ 9\right] $ the
existence of the global classical solutions and the blow-up of the solutions
of the initial boundary value problem and Cauchy problem $\left( 1.3\right)
-\left( 1.4\right) $ are obtained.\ Here, by inspiring $\left[ 8\right] $
and $\left[ 9\right] ,$ the Cauchy problem for Boussinesq operator equation
is considered. Note that, differential operator equations were studied e.g.
in $\left[ \text{10-42, 53}\right] .$ Cauchy problem for abstract hyperbolic
equations were treated e.g. in $\left[ \text{11-20}\right] $ and for
abstract Boussinesq equations studied in $\left[ \text{35, 36}\right] .$ In
this paper, we obtain the local existence and uniqueness of small-amplitude
solution of the Cauchy problem for abstract Boussinesq equations with
general elliptic principal part$.$ The strategy is to express the abstract
Boussinesq equation as an integral equation with operator coefficient. To
treat the nonlinearity as a small perturbation of the linear part of the
equation, the contraction mapping theorem is used. Also, a priori estimates
on $E-$valued $L^{p}$ norms of solutions of the linearized version are
utilized. The key step is the derivation of the uniform estimate of the
solutions of the linearized Boussinesq-operator equation. Modern analysis
methods, particularly abstract harmonic analysis, operator theory,
interpolation of Banach Spaces, embedding theorems in abstract Sobolev-Lions
spaces are the main tools implemented to carry out the analysis.

In order to state our results precisely, we introduce some notations and
some function spaces.

\begin{center}
\textbf{Definitions and} \textbf{Background}
\end{center}

Let $E$ be a Banach space. $L^{p}\left( \Omega ;E\right) $ denotes the space
of strongly measurable $E$-valued functions that are defined on the
measurable subset $\Omega \subset R^{n}$ with the norm

\[
\left\Vert f\right\Vert _{L^{p}}=\left\Vert f\right\Vert _{L^{p}\left(
\Omega ;E\right) }=\left( \int\limits_{\Omega }\left\Vert f\left( x\right)
\right\Vert _{E}^{p}dx\right) ^{\frac{1}{p}},1\leq p<\infty ,\text{ } 
\]

\[
\left\Vert f\right\Vert _{L^{\infty }}\ =\text{ess}\sup\limits_{x\in \Omega
}\left\Vert f\left( x\right) \right\Vert _{E}. 
\]%
The Banach space\ $E$ is called an UMD-space if\ the Hilbert operator%
\[
\left( Hf\right) \left( x\right) =\lim\limits_{\varepsilon \rightarrow
0}\int\limits_{\left\vert x-y\right\vert >\varepsilon }\frac{f\left(
y\right) }{x-y}dy 
\]%
\ is bounded in $L^{p}\left( R,E\right) ,$ $p\in \left( 1,\infty \right) $ (
see. e.g. $\left[ 43\right] $ ). UMD spaces include e.g. $L^{p}$, $l_{p}$
spaces and Lorentz spaces $L_{pq}$ for $p$, $q\in \left( 1,\infty \right) $.

Let \ 
\[
S_{\psi }=\left\{ \lambda \in \mathbb{C}\text{, }\left\vert \arg \lambda
\right\vert \leq \omega ,\text{ }0\leq \omega <\pi \right\} , 
\]%
\[
S_{\omega ,\varkappa }=\left\{ \lambda \in S_{\omega }\text{, }\left\vert
\lambda \right\vert >\varkappa >0\right\} \text{ }. 
\]

A closed linear operator\ $A$ is said to be positive in a Banach\ space $E$
if $D\left( A\right) $ is dense on $E$ and $\left\Vert \left( A+\lambda
I\right) ^{-1}\right\Vert _{B\left( E\right) }\leq M\left( 1+\left\vert
\lambda \right\vert \right) ^{-1}$ for any $\lambda \in S_{\omega },$ $0\leq
\omega <\pi ,$ where $I$ is the identity operator in $E,$ $B\left( E\right) $
is the space of bounded linear operators in $E;$ $D\left( A\right) $ denote
domain of the operator $A.$ It is known $\left[ \text{44, \S 1.15.1}\right] $
that there exist fractional powers\ $A^{\theta }$ of a positive operator $A.$
Let $E\left( A^{\theta }\right) $ denote the space $D\left( A^{\theta
}\right) $ with the graphical norm 
\[
\left\Vert u\right\Vert _{E\left( A^{\theta }\right) }=\left( \left\Vert
u\right\Vert ^{p}+\left\Vert A^{\theta }u\right\Vert ^{p}\right) ^{\frac{1}{p%
}},1\leq p<\infty ,\text{ }0<\theta <\infty . 
\]

A closed linear operator\ $A$ in a Banach\ space $E$ belong to $\sigma
\left( C_{0},\omega ,E\right) $ (see $\left[ 11\right] $, \S\ 11.2 ) if $%
D\left( A\right) $ is dense on $E,$ the resolvent $\left( A-\lambda
^{2}I\right) ^{-1}$ exists for $\func{Re}\lambda >\omega $ and 
\[
\left\Vert \left( A-\lambda ^{2}I\right) ^{-1}\right\Vert _{B\left( E\right)
}\leq C_{0}\left\vert \func{Re}\lambda -\omega \right\vert ^{-1}\text{. } 
\]

Let $E_{1}$ and $E_{2}$ be two Banach spaces. $\left( E_{1},E_{2}\right)
_{\theta ,p}$, $0<\theta <1,1\leq p\leq \infty $ denotes the interpolation
spaces obtained from $\left\{ E_{1},E_{2}\right\} $ by $K$-method \ $\left[ 
\text{44, \S 1.3.2}\right] $.

Let $\mathbb{N}$ denote the set of natural numbers. A set $\Phi \subset
B\left( E_{1},E_{2}\right) $ is called $R$-bounded ( see e.g. $\left[ 10%
\right] $ ) if there is a positive constant $C$ such that for all $%
T_{1},T_{2},...,T_{m}\in \Phi $ and $u_{1,}u_{2},...,u_{m}\in E_{1},$ $m\in 
\mathbb{N}$ 
\[
\int\limits_{\Omega }\left\Vert \sum\limits_{j=1}^{m}r_{j}\left( y\right)
T_{j}u_{j}\right\Vert _{E_{2}}dy\leq C\int\limits_{\Omega }\left\Vert
\sum\limits_{j=1}^{m}r_{j}\left( y\right) u_{j}\right\Vert _{E_{1}}dy, 
\]

where $\left\{ r_{j}\right\} $ is a sequence of independent symmetric $%
\left\{ -1,1\right\} $-valued random variables on $\Omega $. The smallest $C$
for which the above estimate holds is called a $R$-bound of the collection $%
\Phi $ and denoted by $R\left( \Phi \right) .$

Let $h$ be same parameter with $h\in Q\subset \mathbb{C}.$ A set $\Phi
_{h}\subset B\left( E_{1},E_{2}\right) $ is called uniform $R$-bounded if
there is a constant $C$ independent on $h$ such that 
\[
\int\limits_{\Omega }\left\Vert \sum\limits_{j=1}^{m}r_{j}\left( y\right)
T_{j}\left( h\right) u_{j}\right\Vert _{E_{2}}dy\leq C\int\limits_{\Omega
}\left\Vert \sum\limits_{j=1}^{m}r_{j}\left( y\right) u_{j}\right\Vert
_{E_{1}}dy. 
\]

for all $T_{1}\left( h\right) ,T_{2}\left( h\right) ,...,T_{m}\left(
h\right) \in \Phi _{h}$ and $u_{1,}u_{2},...,u_{m}\in E_{1},$ $m\in \mathbb{N%
}$. It is implies that $\sup\limits_{h\in Q}R\left( \Phi _{h}\right) \leq C$.

The positive operator $A$ is said to be $R$-positive in a Banach space $E$
if the set $L_{A}=\left\{ \xi \left( A+\xi \right) ^{-1}\text{: }\xi \in
S_{\omega }\right\} $, $0\leq \omega <\pi $ is $R$-bounded.

Let

\[
\alpha =\left( \alpha _{1},\alpha _{2},...,\alpha _{n}\right) \text{, }%
D^{\alpha }=\frac{\partial ^{\left\vert \alpha \right\vert }}{\partial
x_{1}^{\alpha _{1}}\partial x_{2}^{\alpha _{2}}...\partial x_{n}^{\alpha
_{n}}},\text{ }\left\vert \alpha \right\vert =\sum\limits_{k=1}^{n}\alpha
_{k}\text{.} 
\]

Let $E_{0}$ and $E$ be two Banach spaces and $E_{0}$ is continuously and
densely embedded into $E$. Let $\Omega $ be a domain in $R^{n}$ and $m$ be a
positive integer$.$\ $W^{m,p}\left( \Omega ;E_{0},E\right) $ denotes the
space of all functions $u\in L^{p}\left( \Omega ;E_{0}\right) $ that have
the generalized derivatives $\frac{\partial ^{m}u}{\partial x_{k}^{m}}\in
L^{p}\left( \Omega ;E\right) ,$ $1\leq p\leq \infty $ with the norm 
\[
\ \left\Vert u\right\Vert _{W^{m,p}\left( \Omega ;E_{0},E\right)
}=\left\Vert u\right\Vert _{L^{p}\left( \Omega ;E_{0}\right)
}+\sum\limits_{k=1}^{n}\left\Vert \frac{\partial ^{m}u}{\partial x_{k}^{m}}%
\right\Vert _{L^{p}\left( \Omega ;E\right) }<\infty . 
\]%
\ \ $\ \ $ For $E_{0}=E$ the space $W^{m,p}\left( \Omega ;E_{0},E\right) $
denotes by $W^{m,p}\left( \Omega ;E\right) .$

Let $L^{s,p}\left( R^{n};E\right) $, $-\infty <s<\infty $ denotes the $E-$%
valued Liouville-Sobolev space of order $s$ which is defined as: 
\[
L^{s,p}=L^{s,p}\left( R^{n};E\right) =\left( I-\Delta \right) ^{-\frac{s}{2}%
}L^{p}\left( R^{n};E\right) 
\]%
with the norm 
\[
\left\Vert u\right\Vert _{L^{s,p}}=\left\Vert \left( I-\Delta \right) ^{%
\frac{s}{2}}u\right\Vert _{L^{p}\left( R^{n};E\right) }. 
\]%
It clear that $L^{0,p}\left( R^{n};E\right) =L^{p}\left( R^{n};E\right) .$
It is known that if $E$ is a UMD space, then $L^{m,p}\left( R^{n};E\right)
=W^{m,p}\left( R^{n};E\right) $ for positive integer $m$ ( see e.g. $\left[ 
\text{45, \S\ 15}\right] $ )$.$ $L^{s,p}\left( R^{n};E_{0},E\right) $ denote
the Liouville-Lions type space i.e., 
\[
L^{s,p}\left( R^{n};E_{0},E\right) =\left\{ u\in L^{s,p}\left(
R^{n};E\right) \cap L^{q}\left( R^{n};E_{0}\right) \right. \text{, } 
\]%
\[
\left. \left\Vert u\right\Vert _{L^{s,p}\left( R^{n};E_{0},E\right)
}=\left\Vert u\right\Vert _{L^{p}\left( R^{n};E_{0}\right) }+\left\Vert
u\right\Vert _{L^{s,p}\left( R^{n};E\right) }<\infty \right\} . 
\]%
Let $S\left( R^{n};E\right) $ denote $E$-valued Schwartz class, i.e., the
space of $E$-valued rapidly decreasing smooth functions on $R^{n},$ equipped
with its usual topology generated by seminorms. Let $S^{^{\prime }}\left(
R^{n};E\right) $ denote the space of all\ continuous linear operators $%
L:S\left( R^{n};E\right) \rightarrow E,$ equipped with the bounded
convergence topology. Recall $S\left( R^{n};E\right) $ is norm dense in $%
L^{p}\left( R^{n};E\right) $ when $1\leq p<\infty .$

Let $L_{q}^{\ast }\left( E\right) $ denote the space of all $E-$valued
function space such that 
\[
\left\Vert u\right\Vert _{L_{q}^{\ast }\left( E\right) }=\left(
\int\limits_{0}^{\infty }\left\Vert u\left( t\right) \right\Vert _{E}^{q}%
\frac{dt}{t}\right) ^{\frac{1}{q}}<\infty ,\text{ }1\leq q<\infty ,\text{ }%
\left\Vert u\right\Vert _{L_{\infty }^{\ast }\left( E\right)
}=\sup_{0<t<\infty }\left\Vert u\left( t\right) \right\Vert _{E}. 
\]%
Let $s=\left( s_{1},s_{2},...,s_{n}\right) $ and $s_{k}>0$. Let $F$ denote
the Fourier transform. Fourier-analytic representation of $E-$valued Besov
space on $R^{n}$ are defined as:%
\[
B_{p,q}^{s}\left( R^{n};E\right) =\left\{ u\in S^{^{\prime }}\left(
R^{n};E\right) ,\right. \text{ } 
\]%
\[
\left\Vert u\right\Vert _{B_{p,q}^{s}\left( R^{n};E\right) }=\left\Vert
F^{-1}\sum\limits_{k=1}^{n}t^{\varkappa _{k}-s_{k}}\left( 1+\left\vert \xi
_{k}\right\vert ^{\varkappa _{k}}\right) e^{-t\left\vert \xi \right\vert
^{2}}Fu\right\Vert _{L_{q}^{\ast }\left( L_{p}\left( R^{n};E\right) \right) }%
\text{,} 
\]%
\[
\left. p\in \left( 1,\infty \right) \text{, }q\in \left[ 1,\infty \right] 
\text{, }\varkappa _{k}>s_{k}\right\} . 
\]

It should be note that, the norm of Besov space does not depends on $%
\varkappa _{k}.$ See ( $\left[ \text{44, \S\ 2.3}\right] $ for the scalar
case, i.e., $E=\mathbb{C}$ ).

Let $B_{p,q}^{s}\left( R^{n};E_{0},E\right) $ denote the space $L^{p}\left(
R^{n};E_{0}\right) \cap B_{p,q}^{s}\left( R^{n};E\right) $ with the norm 
\[
\left\Vert u\right\Vert _{B_{p,q}^{s}\left( R^{n};E_{0},E\right)
}=\left\Vert u\right\Vert _{L^{p}\left( R^{n};E_{0}\right) }+\left\Vert
u\right\Vert _{B_{p,q}^{s}\left( R^{n};E\right) }<\infty . 
\]

The embedding theorems in vector valued spaces play a key role in the theory
of DOEs. For estimating lower order derivatives we use following embedding
theorem that is obtained from $\left[ \text{32, Theorem 1}\right] $:

\textbf{Theorem A}$_{1}$. Suppose the following conditions are satisfied:

(1) $E$ is a UMD space and\ $A$ is an $R$-positive operator in $E;$

(2)\ $\alpha =\left( \alpha _{1},\alpha _{2},...,\alpha _{n}\right) $ is a $%
n $-tuples of nonnegative integer number\ and $s$ is a positive number such
that

$\varkappa =\frac{\left\vert \alpha \right\vert +n\left( \frac{1}{p}-\frac{1%
}{q}\right) }{s}\leq 1,$ $0\leq \mu \leq 1-\varkappa $, $1<p\leq q<\infty ;$ 
$0<h\leq h_{0},$ where $h_{0}$ is a fixed positive number;

Then the embedding $D^{\alpha }L^{s,p}\left( R^{n};E\left( A\right)
,E\right) \subset L^{q}\left( R^{n};E\left( A^{1-\varkappa -\mu }\right)
\right) $ is continuous and for $u\in L^{s,p}\left( R^{n};E\left( A\right)
,E\right) $ the following uniform estimate holds 
\[
\left\Vert D^{\alpha }u\right\Vert _{L^{q}\left( R^{n};E\left(
A^{1-\varkappa -\mu }\right) \right) }\leq h^{\mu }\left\Vert u\right\Vert
_{L^{s,p}\left( R^{n};E\left( A\right) ,E\right) }+h^{-\left( 1-\mu \right)
}\left\Vert u\right\Vert _{L^{p}\left( R^{n};E\right) }. 
\]

In a similar way as $\left[ \text{31, Theorem A}_{0}\right] $ and by
reasoning as $\left[ \text{46, Theorem 3.7}\right] $ we obtain:

\textbf{Proposition A}$_{1}.$ Let $1<p\leq q\leq \infty $ and $E$ be $UMD$
space. Suppose $\Psi _{h}\in C^{n}\left( R^{n}\backslash \left\{ 0\right\}
;B\left( E\right) \right) $ and there is a positive constant $K$ such that

\[
\text{ }\sup\limits_{h\in Q}R\left( \left\{ \left\vert \xi \right\vert
^{\left\vert \beta \right\vert +n\left( \frac{1}{p}-\frac{1}{q}\right)
}D^{\beta }\Psi _{h}\left( \xi \right) \text{: }\xi \in R^{n}\backslash
\left\{ 0\right\} ,\text{ }\beta _{k}\in \left\{ 0,1\right\} \right\}
\right) \leq K.\text{ } 
\]

Then $\Psi _{h}$ is a uniformly bounded collection of Fourier multiplier
from $L^{p}\left( R^{n};E\right) $ to $L^{q}\left( R^{n};E\right) .$

\textbf{Proof. }First, in a similar way as in $\left[ \text{31, Theorem A}%
_{0}\right] $ we show that $\Psi _{h}$ is a uniformly bounded collection of
Fourier multiplier from $L^{p}\left( R^{n};E\right) $ to $L^{p}\left(
R^{n};E\right) .$ Moreover, by Theorem A$_{1}$ we get that, for $s\geq
n\left( \frac{1}{p}-\frac{1}{q}\right) $ the embedding $L^{s,p}\left(
R^{n};E\right) \subset L^{q}\left( R^{n};E\right) $ is continuous. From
these two fact we obtain the conclusion.

Sometimes we use one and the same symbol $C$ without distinction in order to
denote positive constants which may differ from each other even in a single
context. When we want to specify the dependence of such a constant on a
parameter, say $\alpha $, we write $C_{\alpha }$.

The paper is organized as follows: In Section 1, some definitions and
background are given. In Section 2, we obtain the existence of unique
solution and a priory estimates for solution of the linearized problem $%
(1.1) $-$\left( 1.2\right) .$ In Section 3, we show the existence and
uniqueness of local strong solution of the problem $(1.1)$-$\left(
1.2\right) $. In Section 4, the existence, uniqueness and a priory estimates
for solution of Cauchy problem for finite and infinite system of Boussinesg
equation is derived.

Sometimes we use one and the same symbol $C$ without distinction in order to
denote positive constants which may differ from each other even in a single
context. When we want to specify the dependence of such a constant on a
parameter, say $h$, we write $C_{h}$.

\begin{center}
\textbf{2. Estimates for linearized equation}
\end{center}

In this section, we make the necessary estimates for solutions of initial
value problems for the linearized abstract Boussinesq equation%
\begin{equation}
u_{tt}-Lu_{tt}+Au=g\left( x,t\right) ,\text{ }x\in R^{n},\text{ }t\in \left(
0,T\right) ,  \tag{2.1}
\end{equation}%
\begin{equation}
u\left( x,0\right) =\varphi \left( x\right) ,\text{ }u_{t}\left( x,0\right)
=\psi \left( x\right) ,  \tag{2.2}
\end{equation}

\bigskip where 
\[
Lu=\dsum\limits_{i,j=1}^{2}a_{ij}\frac{\partial ^{2}u}{\partial
x_{i}\partial x_{j}},\text{ }a_{ij}\in \mathbb{C}. 
\]%
\textbf{Condition 2.0. }Assume $L$ is an elliptic operator, i.e, there are
positive constants $M_{1}$ and $M_{2}$ such that $\ M_{1}\left\vert \xi
\right\vert ^{2}\leq L\left( \xi \right) \leq M_{2}\left\vert \xi
\right\vert ^{2}$ for $\xi =\left( \xi _{1},\xi _{2},...\xi _{n}\right) \in
R^{n}$, where%
\[
\left\vert \xi \right\vert ^{2}=\dsum\limits_{k=1}^{n}\xi _{k}^{2}\text{, }%
L\left( \xi \right) =\dsum\limits_{i,j=1}^{2}a_{ij}\xi _{i}\xi _{j}. 
\]%
Let 
\[
X_{p}=L^{p}\left( R^{n};E\right) \text{, }Y^{s,p}=L^{s,p}\left(
R^{n};E\right) ,\text{ }Y_{1}^{s,p}= 
\]

\[
L^{s,p}\left( R^{n};E\right) \cap L^{1}\left( R^{n};E\right) \text{, }%
Y_{\infty }^{s,p}=L^{s,p}\left( R^{n};E\right) \cap L^{\infty }\left(
R^{n};E\right) . 
\]%
\textbf{Condition 2.1. }Assume:

(1) $E$ is an UMD space and linear operator $A$ belongs to $\sigma \left(
C_{0},\omega ,E\right) $;

(2) $\varphi ,$ $\psi \in Y_{\infty }^{s,p}$ and $g\left( .,t\right) $ $\in
Y_{\infty }^{s,p}$ for $t\in \left( 0,T\right) $ and $s>\frac{n}{p}$ for $%
1<p<\infty $.

First we need the following lemmas

\textbf{Lemma 2.1. }Suppose the Conditions 2.0, 2.1 hold.Then problem $%
\left( 2.1\right) -\left( 2.2\right) $ has a generalized solution.

\textbf{Proof. }By using of Fourier transform we get from $(2.1)-(2.2)$:%
\begin{equation}
\hat{u}_{tt}\left( \xi ,t\right) +A_{\xi }\hat{u}\left( \xi ,t\right) =\left[
1+L\left( \xi \right) \right] ^{-1}\hat{g}\left( \xi ,t\right) ,\text{ } 
\tag{2.3}
\end{equation}%
\[
\hat{u}\left( \xi ,0\right) =\hat{\varphi}\left( \xi \right) ,\text{ }\hat{u}%
_{t}\left( \xi ,0\right) =\hat{\psi}\left( \xi \right) ,\text{ }\xi \in
R^{n},\text{ }t\in \left( 0,T\right) , 
\]

where $\hat{u}\left( \xi ,t\right) $ is a Fourier transform of $u\left(
x,t\right) $ with respect to $x,$ where%
\[
A_{\xi }=\left[ 1+L\left( \xi \right) \right] ^{-1}A\text{, }\xi \in R^{n}. 
\]

By virtue of $\left[ \text{11, \S 11.2, 11.4}\right] $ ( or $\left[ \text{%
12-20}\right] $ ) we obtain that $A_{\xi }$ is a generator of a strongly
continuous cosine operator function\ and problem $(2.3)$ has a unique
solution for all $\xi \in R^{n},$ moreover, the solution can be written as%
\begin{equation}
\hat{u}\left( \xi ,t\right) =C\left( t,\xi ,A\right) \hat{\varphi}\left( \xi
\right) +S\left( t,\xi ,A\right) \hat{\psi}\left( \xi \right) +  \tag{2.4}
\end{equation}%
\[
\dint\limits_{0}^{t}S\left( t-\tau ,\xi ,A\right) \left[ 1+L\left( \xi
\right) \right] ^{-1}\hat{g}\left( \xi ,\tau \right) d\tau ,\text{ }t\in
\left( 0,T\right) , 
\]%
where $C\left( t,\xi ,A\right) $ is a cosine and $S\left( t,\xi ,A\right) $
is a sine operator-functions (see e.g. $\left[ 11\right] $) generated by
parameter dependent operator $A_{\xi }.$ From $\left( 2.4\right) $ we get
that, the solution of the problem $(2.1)-(2.2)$ can be expressed as 
\[
u\left( x,t\right) =S_{1}\left( t,A\right) \varphi \left( x\right)
+S_{2}\left( t,A\right) \psi \left( x\right) + 
\]%
\begin{equation}
+\left( 2\pi \right) ^{-\frac{1}{n}}\dint\limits_{R^{n}}\dint%
\limits_{0}^{t}e^{ix\xi }S\left( t-\tau ,\xi ,A\right) \left[ 1+L\left( \xi
\right) \right] ^{-1}\hat{g}\left( \xi ,\tau \right) d\tau d\xi ,\text{ }%
t\in \left( 0,T\right) ,  \tag{2.5}
\end{equation}%
where $S_{1}\left( t,A\right) $ and $S_{2}\left( t,A\right) $ are linear
operators in $E$ defined by 
\begin{equation}
S_{1}\left( t,A\right) \varphi =\left( 2\pi \right) ^{-\frac{1}{n}%
}\dint\limits_{R^{n}}e^{ix\xi }C\left( t,\xi ,A\right) \hat{\varphi}\left(
\xi \right) d\xi ,\text{ }  \tag{2.6}
\end{equation}%
\[
S_{2}\left( t,A\right) \psi =\left( 2\pi \right) ^{-\frac{1}{n}%
}\dint\limits_{R^{n}}e^{ix\xi }S\left( t,\xi ,A\right) \hat{\psi}\left( \xi
\right) d\xi .\text{ } 
\]

\textbf{Lemma 2.2. }Suppose the Conditions 2.0, 2.1 hold. Then the solution
of the problem $\left( 2.1\right) -\left( 2.2\right) $ satisfies the
following estimate 
\begin{equation}
\left( \left\Vert u\right\Vert _{X_{\infty }}+\left\Vert u_{t}\right\Vert
_{X_{\infty }}\right) \leq C\left( \left\Vert \varphi \right\Vert
_{Y^{s,p}}+\left\Vert \varphi \right\Vert _{X_{1}}\right.  \tag{2.7}
\end{equation}

\[
+\left\Vert \psi \right\Vert _{Y^{s,p}}+\left\Vert \psi \right\Vert
_{X_{1}}+\dint\limits_{0}^{t}\left( \left\Vert g\left( .,\tau \right)
\right\Vert _{Y^{s,p}}+\left\Vert g\left( .,\tau \right) \right\Vert
_{X_{1}}\right) d\tau 
\]%
uniformly in $t\in \left[ 0,T\right] .$

\textbf{Proof. }Let $N\in \mathbb{N}$ and 
\[
\Pi _{N}=\left\{ \xi :\xi \in R^{n},\text{ }\left\vert \xi \right\vert \leq
N\right\} ,\text{ }\Pi _{N}^{\prime }=\left\{ \xi :\xi \in R^{n},\text{ }%
\left\vert \xi \right\vert \geq N\right\} . 
\]%
It is clear to see that

\[
\left\Vert u\left( .,t\right) \right\Vert _{L^{\infty }\left( R^{n};E\right)
}=\left\Vert F^{-1}\hat{u}\left( \xi ,t\right) \right\Vert _{L^{\infty
}\left( R^{n};E\right) }\leq 
\]

\[
\left\Vert F^{-1}C\left( t,\xi ,A\right) \hat{\varphi}\left( \xi \right)
\right\Vert _{X_{\infty }}+\left\Vert F^{-1}S\left( t,\xi ,A\right) \hat{\psi%
}\left( \xi \right) \right\Vert _{X_{\infty }}\leq 
\]%
\begin{equation}
\left\Vert F^{-1}C\left( t,\xi ,A\right) \hat{\varphi}\left( \xi \right)
\right\Vert _{L^{\infty }\left( \Pi _{N};E\right) }+\left\Vert F^{-1}S\left(
t,\xi ,A\right) \hat{\psi}\left( \xi \right) \right\Vert _{L^{\infty }\left(
\Pi _{N};E\right) }+  \tag{2.8}
\end{equation}%
\[
\left\Vert F^{-1}C\left( t,\xi ,A\right) \hat{\varphi}\left( \xi \right)
\right\Vert _{L^{\infty }\left( \Pi _{N}^{\prime };E\right) }+\left\Vert
F^{-1}S\left( t,\xi ,A\right) \hat{\psi}\left( \xi \right) \right\Vert
_{L^{\infty }\left( \Pi _{N}^{\prime };E\right) }, 
\]%
\[
\left\Vert F^{-1}C\left( t,\xi ,A\right) \hat{\varphi}\left( \xi \right)
\right\Vert _{L^{\infty }\left( \Pi _{N}^{\prime };E\right) }+\left\Vert
F^{-1}S\left( t,\xi ,A\right) \hat{\psi}\left( \xi \right) \right\Vert
_{L^{\infty }\left( \Pi _{N}^{\prime };E\right) }= 
\]%
\begin{equation}
=\left\Vert F^{-1}\left( 1+L\left( \xi \right) \right) ^{-\frac{s}{2}%
}C\left( t,\xi ,A\right) \left( 1+L\left( \xi \right) \right) ^{\frac{s}{2}}%
\hat{\varphi}\left( \xi \right) \right\Vert _{L^{\infty }\left( \Pi
_{N}^{\prime };E\right) }+  \tag{2.9}
\end{equation}%
\[
\left\Vert F^{-1}\left( 1+L\left( \xi \right) \right) ^{-\frac{s}{2}}S\left(
t,\xi ,A\right) \left( 1+L\left( \xi \right) \right) ^{\frac{s}{2}}\hat{\psi}%
\left( \xi \right) \right\Vert _{L^{\infty }\left( \Pi _{N}^{\prime
};E\right) }. 
\]%
Using the H\"{o}lder inequality we have 
\begin{equation}
\left\Vert F^{-1}C\left( t,\xi ,A\right) \hat{\varphi}\left( \xi \right)
\right\Vert _{L^{\infty }\left( \Pi _{N};E\right) }+\left\Vert F^{-1}S\left(
t,\xi ,A\right) \hat{\psi}\left( \xi \right) \right\Vert _{L^{\infty }\left(
\Pi _{N};E\right) }\leq  \tag{2.10}
\end{equation}

\[
C\left[ \left\Vert \varphi \right\Vert _{X_{1}}+\left\Vert \psi \right\Vert
_{X_{1}}\right] . 
\]%
By using the resolvent properties of operator $A$, representation of $%
C\left( t,\xi ,A\right) $, $S\left( t,\xi ,A\right) $ and\ the Conditon 2.0
we get 
\[
\left\vert \xi \right\vert ^{\left\vert \alpha \right\vert +\frac{n}{p}%
}\left\Vert D^{\alpha }\left[ \left( 1+L\left( \xi \right) \right) ^{-\frac{s%
}{2}}C\left( t,\xi ,A\right) \right] \right\Vert _{B\left( E\right) }\leq
C_{1}, 
\]%
\ 
\begin{equation}
\left\vert \xi \right\vert ^{\left\vert \alpha \right\vert +\frac{n}{p}%
}\left\Vert D^{\alpha }\left[ \left( 1+L\left( \xi \right) \right) ^{-\frac{s%
}{2}}S\left( t,\xi ,A\right) \right] \right\Vert _{B\left( E\right) }\leq
C_{2},  \tag{2.11}
\end{equation}%
for $s>\frac{n}{p}$\ and all $\alpha =\left( \alpha _{1},\alpha
_{2},...,\alpha _{n}\right) $, $\alpha _{k}\in \left\{ 0,1\right\} $, $\xi
\in R^{n}$, $\xi \neq 0$, $t\in \left[ 0,T\right] .$ By Proposition A$_{1}$
from $\left( 2.9\right) $ and Conditon 2.0 we get that, the operator-valued
functions $\left( 1+L\left( \xi \right) \right) ^{-\frac{s}{2}}C\left( t,\xi
,A\right) ,$ $\left( 1+L\left( \xi \right) \right) ^{-\frac{s}{2}}S\left(
t,\xi ,A\right) $ are $L^{p}\left( R^{n};E\right) \rightarrow L^{\infty
}\left( R^{n};E\right) $ Fourier multipliers uniformly in $t\in \left[ 0,T%
\right] $. Then by Minkowski's inequality for integrals, the semigroups
estimates (see e.g. $\left[ \text{11-12}\right] $ ) and $\left( 2.9\right) $
we obtain%
\begin{equation}
\left\Vert F^{-1}C\left( t,\xi ,A\right) \hat{\varphi}\left( \xi \right)
\right\Vert _{L^{\infty }\left( \Pi _{N}^{\prime };E\right) }+\left\Vert
F^{-1}S\left( t,\xi ,A\right) \hat{\psi}\left( \xi \right) \right\Vert
_{L^{\infty }\left( \Pi _{N}^{\prime };E\right) }\leq  \tag{2.12}
\end{equation}

\[
C\left[ \left\Vert \varphi \right\Vert _{Y^{s,p}}+\left\Vert \psi
\right\Vert _{Y^{s,p}}\right] . 
\]

By reasoning as the above we get 
\begin{equation}
\left\Vert F^{-1}\dint\limits_{0}^{t}S\left( t-\tau ,\xi ,A\right) \left[
1+L\left( \xi \right) \right] ^{-1}\hat{g}\left( \xi ,\tau \right) d\tau
\right\Vert _{X_{\infty }}\leq  \tag{2.13}
\end{equation}

\[
C\dint\limits_{0}^{t}\left( \left\Vert g\left( .,\tau \right) \right\Vert
_{Y^{s}}+\left\Vert g\left( .,\tau \right) \right\Vert _{X_{1}}\right) d\tau
. 
\]%
By differentiating, in view of $\left( 2.6\right) $ we obtain from $\left(
2.5\right) $ the estimate of type $\left( 2.10\right) ,$ $\left( 2.12\right)
,$ $\left( 2.13\right) $ for $u_{t}.$

\bigskip Then by using $\left( 2.10\right) ,$ $\left( 2.12\right) ,$ $\left(
2.13\right) $ we get the estimate $\left( 2.7\right) .$

\textbf{Lemma 2.3. }Assume the Conditions 2.0, 2.1 are hold. Then the
solution of the problem $\left( 2.1\right) -\left( 2.2\right) $ satisfies
the following uniform estimate%
\begin{equation}
\left( \left\Vert u\right\Vert _{Y^{s,p}}+\left\Vert u_{t}\right\Vert
_{Y^{s,p}}\right) \leq C\left( \left\Vert \varphi \right\Vert
_{Y^{s,p}}+\left\Vert \psi \right\Vert
_{Y^{s,p}}+\dint\limits_{0}^{t}\left\Vert g\left( .,\tau \right) \right\Vert
_{Y^{s,p}}d\tau \right) .  \tag{2.14}
\end{equation}

\textbf{Proof. }From $\left( 2.4\right) $ we have the following estimate 
\[
\left( \left\Vert F^{-1}\left( 1+\left\vert \xi \right\vert ^{2}\right) ^{%
\frac{s}{2}}\hat{u}\right\Vert _{X_{p}}+\left\Vert F^{-1}\left( 1+\left\vert
\xi \right\vert ^{2}\right) ^{\frac{s}{2}}\hat{u}_{t}\right\Vert
_{X_{p}}\right) \leq 
\]

\begin{equation}
C\left\{ \left\Vert F^{-1}\left( 1+\left\vert \xi \right\vert \right) ^{%
\frac{s}{2}}C\left( t,\xi ,A\right) \hat{\varphi}\right\Vert _{X_{p}}\right.
+\left\Vert F^{-1}\left( 1+\left\vert \xi \right\vert \right) ^{\frac{s}{2}%
}S\left( t,\xi ,A\right) \hat{\psi}\right\Vert _{X_{p}}  \tag{2.15}
\end{equation}

\[
+\left. \dint\limits_{0}^{t}\left\Vert \left( 1+\left\vert \xi \right\vert
\right) ^{\frac{s}{2}}S\left( t-\tau ,\xi ,A\right) \left[ 1+L\left( \xi
\right) \right] ^{-1}g\left( .,\tau \right) \right\Vert _{X_{p}}d\tau
\right\} . 
\]

\bigskip By construction of operator-valued functions $C\left( t,\xi
,A\right) $, $S\left( t,\xi ,A\right) $ and in view of Proposition A$_{1}$
and Conditon 2.0\ we get that $C\left( t,\xi ,A\right) $ and $S\left( t,\xi
,A\right) $ are $L^{p}\left( R^{n};E\right) $ Fourier multipliers uniformly
in $t\in \left[ 0,T\right] .$ So, the estimate $\left( 2.15\right) $ by
using the Minkowski's inequality for integrals implies $\left( 2.14\right) .$

From Lemmas 2.1-2.3 we obtain

\textbf{Theorem 2.1. }Let the Condition 2.1 hold. Then problem $\left(
2.1\right) -\left( 2.2\right) $ has a unique solution $u\in C^{\left(
2\right) }\left( \left[ 0,T\right] ;Y_{1}^{s,p}\right) $\textbf{\ }and the
following estimates holds

\begin{equation}
\left\Vert u\right\Vert _{X_{\infty }}+\left\Vert u_{t}\right\Vert
_{X_{\infty }}\leq C\left( \left\Vert \varphi \right\Vert
_{Y^{s,p}}+\left\Vert \varphi \right\Vert _{X_{1}}\right.  \tag{2.16}
\end{equation}

\[
+\left\Vert \psi \right\Vert _{Y^{s,p}}+\left\Vert \psi \right\Vert
_{X_{1}}+\dint\limits_{0}^{t}\left( \left\Vert g\left( .,\tau \right)
\right\Vert _{Y^{s,p}}+\left\Vert g\left( .,\tau \right) \right\Vert
_{X_{1}}\right) d\tau , 
\]

\begin{equation}
\left\Vert u\right\Vert _{Y^{s,p}}+\left\Vert u_{t}\right\Vert
_{Y^{s,p}}\leq C\left( \left\Vert \varphi \right\Vert _{Y^{s,p}}+\left\Vert
\psi \right\Vert _{Y^{s,p}}+\dint\limits_{0}^{t}\left\Vert g\left( .,\tau
\right) \right\Vert _{Y^{s,p}}d\tau \right)  \tag{2.17}
\end{equation}%
uniformly in $t\in \left[ 0,T\right] .$

\textbf{Proof. }From Lemma 2.1 we obtain that, problem $\left( 2.1\right)
-\left( 2.2\right) $ has a unique generalized solution. From the
representation of solution $\left( 2.5\right) $\ and Lemmas 2.2, 2.3 we get
that there is a solution $u\in C^{\left( 2\right) }\left( \left[ 0,T\right]
;Y_{1}^{s}\right) $ and estimates $\left( 2.16\right) $, $\left( 2.17\right) 
$ hold.

\begin{center}
\textbf{3. Initial value problem for nonlinear equation}
\end{center}

In this section, we will show the local existence and uniqueness of solution
for the Cauchy problem $(1.1),(1.2).$

For the study of the nonlinear problem $\left( 1.1\right) -\left( 1.2\right) 
$ we need the following lemmas

\textbf{Lemma 3.1} (Abstract Nirenberg's inequality). Let $E$ be an $UMD$
space. Assume that $u\in L_{p}\left( \Omega ;E\right) $, $D^{m}u$ $\in
L_{q}\left( \Omega ;E\right) $, $p,q\in \left( 1,\infty \right) $. Then for $%
i$ with $0\leq i\leq m,$ $m>\frac{n}{q}$ we have 
\begin{equation}
\left\Vert D^{i}u\right\Vert _{r}\leq C\left\Vert u\right\Vert _{p}^{1-\mu
}\dsum\limits_{k=1}^{n}\left\Vert D_{k}^{m}u\right\Vert _{q}^{\mu }, 
\tag{3.1}
\end{equation}%
where%
\[
\frac{1}{r}=\frac{i}{m}+\mu \left( \frac{1}{q}-\frac{m}{n}\right) +\left(
1-\mu \right) \frac{1}{p},\text{ }\frac{i}{m}\leq \mu \leq 1. 
\]

\textbf{Proof. }By virtue of interpolation of Banach spaces $\left[ \text{%
44, \S 1.3.2}\right] ,$ in order to prove $\left( 3.1\right) $ for any given 
$i,$ one has only to prove it for the extreme values $\mu =\frac{i}{m}$ and $%
\mu =1$. For the case of $\mu =1$, i.e., $\frac{1}{r}=\frac{i}{m}+\frac{1}{q}%
-\frac{m}{n}$ the estimate $\left( 3.1\right) $ is obtained from Theorem A$%
_{1}$. The case $\mu =\frac{i}{m}$ is derived by reasoning as in $\left[ 
\text{47, \S\ 2 }\right] $ and in replacing absolute value of complex-valued
function $u$ by the $E-$norm of $E$-valued function.

\bigskip Note that, for $E=\mathbb{C}$ the lemma considered by L. Nirenberg $%
\left[ 47\right] .$

Using the chain rule of the composite function, from Lemma 3.1 we can prove
the following result

\textbf{Lemma 3.2. }Let $E$ be an $UMD$ space. Assume that $u\in $ $%
W^{m,p}\left( \Omega ;E\right) \cap L^{\infty }\left( \Omega ;E\right) $,
and $f\left( u\right) $ possesses continuous derivatives up to order $m\geq
1 $. Then $f\left( u\right) -f\left( 0\right) \in W^{m,p}\left( \Omega
;E\right) $ and 
\[
\left\Vert f\left( u\right) -f\left( 0\right) \right\Vert _{p}\leq
\left\Vert f^{^{\left( 1\right) }}\left( u\right) \right\Vert _{\infty
}\left\Vert u\right\Vert _{p}, 
\]

\begin{equation}
\left\Vert D^{k}f\left( u\right) \right\Vert _{p}\leq
C_{0}\dsum\limits_{j=1}^{k}\left\Vert f^{\left( j\right) }\left( u\right)
\right\Vert _{\infty }\left\Vert u\right\Vert _{\infty }^{j-1}\left\Vert
D^{k}u\right\Vert _{p}\text{, }1\leq k\leq m,  \tag{3.2}
\end{equation}%
where $C_{0}$ $\geq 1$ is a constant.

For $E=\mathbb{C}$ the lemma coincide with the corresponding inequality in $%
\left[ 48\right] .$ Let%
\[
\text{ }X=L^{p}\left( R^{n};E\right) ,\text{ }Y=W^{2,p}\left( R^{n};E\left(
A\right) ,E\right) ,\text{ }E_{0}=\left( X,Y\right) _{\frac{1}{2p},p}. 
\]

\bigskip \textbf{Remark 3.1. }By using J.Lions-I. Petree result ( see e.g $%
\left[ \text{49}\right] $ or $\left[ \text{44, \S\ 1.8.}\right] $ ) we
obtain that the map $u\rightarrow u\left( t_{0}\right) $, $t_{0}\in \left[
0,T\right] $ is continuous from $W^{2,p}\left( 0,T;X,Y\right) $ onto $E_{0}$
and there is a constant $C_{1}$ such that 
\[
\left\Vert u\left( t_{0}\right) \right\Vert _{E_{0}}\leq C_{1}\left\Vert
u\right\Vert _{W^{2,p}\left( 0,T;X,Y\right) },\text{ }1\leq p\leq \infty 
\text{.} 
\]

Here, we define the space $Y\left( T\right) =C\left( \left[ 0,T\right]
;Y_{\infty }^{2,p}\right) $ equipped with the norm defined by%
\[
\left\Vert u\right\Vert _{Y\left( T\right) }=\max\limits_{t\in \left[ 0,T%
\right] }\left\Vert u\right\Vert _{Y^{2,p}}+\max\limits_{t\in \left[ 0,T%
\right] }\left\Vert u\right\Vert _{X_{\infty }},\text{ }u\in Y\left(
T\right) . 
\]

It is easy to see that $Y\left( T\right) $ is a Banach space. For $\varphi $%
, $\psi \in Y^{2,p}$, let 
\[
M=\left\Vert \varphi \right\Vert _{Y^{2,p}}+\left\Vert \varphi \right\Vert
_{X_{\infty }}+\left\Vert \psi \right\Vert _{Y^{2,p}}+\left\Vert \psi
\right\Vert _{X_{\infty }}. 
\]

\textbf{Definition 3.1. }For any $T>0$ if $\upsilon ,$ $\psi \in Y_{\infty
}^{2,p}$ and $u$ $\in C\left( \left[ 0,T\right] ;Y_{\infty }^{2,p}\right) $
satisfies the equation $(1.1)-(1.2)$ then $u\left( x,t\right) $ is called
the continuous solution\ or the strong solution of the problem $(1.1)-(1.2).$
If $T<\infty $, then $u\left( x,t\right) $ is called the local strong
solution of the problem $(1.1)-(1.2).$ If $T=\infty $, then $u\left(
x,t\right) $ is called the global strong solution of the problem $%
(1.1)-(1.2) $.

\textbf{Condition 3.1. }Assume:

(1) the operator $A$ generates continuous cosine operator function in UMD
space $E$;

(2) $\varphi ,$ $\psi $ $\in Y_{\infty }^{2,p}$ and $1<p<\infty $ for$\frac{n%
}{p}<2$;

(3) the function $u\rightarrow $ $f\left( u\right) $: $R^{n}\times \left[ 0,T%
\right] \times E_{0}\rightarrow E$ is a measurable in $\left( x,t\right) \in
R^{n}\times \left[ 0,T\right] $ for $u\in E_{0};$ $f\left( x,t,.,.\right) $
is continuous in $u\in E_{0}$ for $x\in R^{n},$ $t\in \left[ 0,T\right] $
and $f\left( u\right) \in C^{\left( 3\right) }\left( E_{0};E\right) .$

\bigskip Main aim of this section is to prove the following result:

\bigskip \textbf{Theorem 3.1. }Let the Condition 3.1 hold. Then problem $%
\left( 1.1\right) -\left( 2.2\right) $ has a unique local strange solution $%
u\in C^{\left( 2\right) }\left( \left[ 0\right. ,\left. T_{0}\right)
;Y_{\infty }^{2,p}\right) $, where $T_{0}$ is a maximal time interval that
is appropriately small relative to $M$. Moreover, if

\begin{equation}
\sup_{t\in \left[ 0\right. ,\left. T_{0}\right) }\left( \left\Vert
u\right\Vert _{Y^{2,p}}+\left\Vert u\right\Vert _{X_{\infty }}+\left\Vert
u_{t}\right\Vert _{Y^{2,p}}+\left\Vert u_{t}\right\Vert _{X_{\infty
}}\right) <\infty  \tag{3.3}
\end{equation}%
then $T_{0}=\infty .$

\textbf{Proof. }First, we are going to prove the existence and the
uniqueness of the local continuous solution of the problem $(1.1)-\left(
1.2\right) $ by contraction mapping principle. Consider a map $G$ on $%
Y\left( T\right) $ such that $G(u)$ is the solution of the Cauchy problem%
\begin{equation}
G_{tt}\left( u\right) -LG_{tt}\left( u\right) +AG\left( u\right) =f\left(
G\left( u\right) \right) ,\text{ }x\in R^{n},\text{ }t\in \left( 0,T\right) ,
\tag{3.4}
\end{equation}%
\[
G\left( u\right) \left( x,0\right) =\varphi \left( x\right) ,\text{ }%
G_{t}\left( u\right) \left( x,0\right) =\psi \left( x\right) . 
\]%
From Lemma 3.2 we know that $f(u)\in $ $L^{p}\left( 0,T;Y_{\infty
}^{2,p}\right) $ for any $T>0$. Thus, by Theorem 2.1, problem $\left(
3.4\right) $ has a unique solution which can be written as%
\[
G\left( u\right) \left( t,x\right) =S_{1}\left( t,A\right) \varphi \left(
x\right) +S_{2}\left( t,A\right) \psi \left( x\right) + 
\]%
\begin{equation}
+\dint\limits_{0}^{t}F^{-1}S\left( t-\tau ,\xi ,A\right) \left[ 1+L\left(
\xi \right) \right] ^{-1}\hat{f}\left( u\right) \left( \xi ,\tau \right)
d\tau ,\text{ }t\in \left( 0,T\right) .  \tag{3.5}
\end{equation}

From Lemma 3.2 it is easy to see that the map $G$ is well defined for $f\in
C^{\left( 2\right) }\left( X_{0};E\right) $. We put 
\[
Q\left( M;T\right) =\left\{ u\mid u\in Y\left( T\right) \text{, }\left\Vert
u\right\Vert _{Y\left( T\right) }\leq M+1\right\} . 
\]

First, by reasoning as in $\left[ 9\right] $\ let us prove that the map $G$
has a unique fixed point in $Q\left( M;T\right) .$ For this aim, it is
sufficient to show that the operator $G$ maps $Q\left( M;T\right) $ into $%
Q\left( M;T\right) $ and $G:$ $Q\left( M;T\right) $ $\rightarrow $ $Q\left(
M;T\right) $ is strictly contractive if $T$ is appropriately small relative
to $M.$ Consider the function \ $\bar{f}\left( \xi \right) $: $\left[
0,\right. $ $\left. \infty \right) \rightarrow \left[ 0,\right. $ $\left.
\infty \right) $ defined by 
\[
\ \bar{f}\left( \xi \right) =\max\limits_{\left\vert x\right\vert \leq \xi
}\left\{ \left\Vert f^{\left( 1\right) }\left( x\right) \right\Vert
_{E},\left\Vert f^{\left( 2\right) }\left( x\right) \right\Vert _{E}\text{ }%
\right\} ,\text{ }\xi \geq 0. 
\]

It is clear to see that the function $\bar{f}\left( \xi \right) $ is
continuous and nondecreasing on $\left[ 0,\right. $ $\left. \infty \right) .$
From Lemma 3.2 we have\qquad

\[
\left\Vert f\left( u\right) \right\Vert _{Y^{2,p}}\leq \left\Vert f^{\left(
1\right) }\left( u\right) \right\Vert _{X_{\infty }}\left\Vert u\right\Vert
_{X_{p}}+\left\Vert f^{\left( 1\right) }\left( u\right) \right\Vert
_{X_{\infty }}\left\Vert Du\right\Vert _{X_{p}} 
\]

\begin{equation}
+C_{0}\left[ \left\Vert f^{\left( 1\right) }\left( u\right) \right\Vert
_{X_{\infty }}\left\Vert u\right\Vert _{X_{p}}+\left\Vert f^{\left( 2\right)
}\left( u\right) \right\Vert _{X_{\infty }}\left\Vert u\right\Vert
_{X_{\infty }}\left\Vert D^{2}u\right\Vert _{X_{p}}\right]  \tag{3.6}
\end{equation}

\ 
\[
\leq 2C_{0}\bar{f}\left( M+1\right) \left( M+1\right) \left\Vert
u\right\Vert _{Y^{2,p}}. 
\]

\bigskip By using Theorem 2.1 we obtain from $\left( 3.5\right) $:%
\begin{equation}
\left\Vert G\left( u\right) \right\Vert _{X_{\infty }}\leq \left\Vert
\varphi \right\Vert _{X_{\infty }}+\left\Vert \psi \right\Vert _{X_{\infty
}}+\dint\limits_{0}^{t}\left\Vert f\left( u\left( \tau \right) \right)
\right\Vert _{X_{\infty }},  \tag{3.7}
\end{equation}%
\begin{equation}
\left\Vert G\left( u\right) \right\Vert _{Y^{2,p}}\leq \left\Vert \varphi
\right\Vert _{Y^{2,p}}+\left\Vert \psi \right\Vert
_{Y^{2,p}}+\dint\limits_{0}^{t}\left\Vert f\left( u\left( \tau \right)
\right) \right\Vert _{Y^{2,p}}d\tau .  \tag{3.8}
\end{equation}

Thus, from $\left( 3.6\right) -\left( 3.8\right) $ and Lemma 3.2 we get 
\[
\left\Vert G\left( u\right) \right\Vert _{Y\left( T\right) }\leq M+T\left(
M+1\right) \left[ 1+2C_{0}\left( M+1\right) \bar{f}\left( M+1\right) \right]
. 
\]

If $T$ satisfies 
\begin{equation}
T\leq \left\{ \left( M+1\right) \left[ 1+2C_{0}\left( M+1\right) \bar{f}%
\left( M+1\right) \right] \right\} ^{-1}.  \tag{3.9}
\end{equation}

Then 
\[
\left\Vert Gu\right\Vert _{Y\left( T\right) }\leq M+1. 
\]%
Therefore, if $\left( 3.9\right) $ holds, then $G$ maps $Q\left( M;T\right) $
into $Q\left( M;T\right) .$ Now, we are going to prove that the map $G$ is
strictly contractive. Assume $T>0$ and $u_{1},$ $u_{2}\in $ $Q\left(
M;T\right) $ given. We get%
\[
G\left( u_{1}\right) -G\left( u_{2}\right)
=\dint\limits_{0}^{t}F^{-1}S\left( t-\tau ,\xi ,A\right) \left[ \hat{f}%
\left( u_{1}\right) \left( \xi ,\tau \right) -\hat{f}\left( u_{2}\right)
\left( \xi ,\tau \right) \right] d\tau ,\text{ }t\in \left( 0,T\right) . 
\]

By using the mean value theorem, we obtain%
\[
\hat{f}\left( u_{1}\right) -\hat{f}\left( u_{2}\right) =\hat{f}^{\left(
1\right) }\left( u_{2}+\eta _{1}\left( u_{1}-u_{2}\right) \right) \left(
u_{1}-u_{2}\right) ,\text{ } 
\]

\[
D_{\xi }\left[ \hat{f}\left( u_{1}\right) -\hat{f}\left( u_{2}\right) \right]
=\hat{f}^{\left( 2\right) }\left( u_{2}+\eta _{2}\left( u_{1}-u_{2}\right)
\right) \left( u_{1}-u_{2}\right) D_{\xi }u_{1}+\text{ } 
\]%
\[
\hat{f}^{\left( 1\right) }\left( u_{2}\right) \left( D_{\xi }u_{1}-D_{\xi
}u_{2}\right) , 
\]%
\[
D_{\xi }^{2}\left[ \hat{f}\left( u_{1}\right) -\hat{f}\left( u_{2}\right) %
\right] =\hat{f}^{\left( 3\right) }\left( u_{2}+\eta _{3}\left(
u_{1}-u_{2}\right) \right) \left( u_{1}-u_{2}\right) \left( D_{\xi
}u_{1}\right) ^{2}+\text{ } 
\]%
\[
\hat{f}^{\left( 2\right) }\left( u_{2}\right) \left( D_{\xi }u_{1}-D_{\xi
}u_{2}\right) \left( D_{\xi }u_{1}+D_{\xi }u_{2}\right) + 
\]%
\[
\hat{f}^{\left( 2\right) }\left( u_{2}+\eta _{4}\left( u_{1}-u_{2}\right)
\right) \left( u_{1}-u_{2}\right) D_{\xi }^{2}u_{1}+\hat{f}^{\left( 1\right)
}\left( u_{2}\right) \left( D_{\xi }^{2}u_{1}-D_{\xi }^{2}u_{2}\right) , 
\]%
where $0<\eta _{i}<1,$ $i=1,2,3,4.$ Thus using H\"{o}lder's and Nirenberg's
inequality, we have%
\begin{equation}
\left\Vert \hat{f}\left( u_{1}\right) -\hat{f}\left( u_{2}\right)
\right\Vert _{X_{\infty }}\leq \bar{f}\left( M+1\right) \left\Vert
u_{1}-u_{2}\right\Vert _{X_{\infty }},  \tag{3.10}
\end{equation}%
\begin{equation}
\left\Vert \hat{f}\left( u_{1}\right) -\hat{f}\left( u_{2}\right)
\right\Vert _{X_{p}}\leq \bar{f}\left( M+1\right) \left\Vert
u_{1}-u_{2}\right\Vert _{X_{p}},  \tag{3.11}
\end{equation}%
\begin{equation}
\left\Vert D_{\xi }\left[ \hat{f}\left( u_{1}\right) -\hat{f}\left(
u_{2}\right) \right] \right\Vert _{X_{p}}\leq \left( M+1\right) \bar{f}%
\left( M+1\right) \left\Vert u_{1}-u_{2}\right\Vert _{X_{\infty }}+ 
\tag{3.12}
\end{equation}%
\[
\bar{f}\left( M+1\right) \left\Vert \hat{f}\left( u_{1}\right) -\hat{f}%
\left( u_{2}\right) \right\Vert _{X_{p}}, 
\]%
\[
\left\Vert D_{\xi }^{2}\left[ \hat{f}\left( u_{1}\right) -\hat{f}\left(
u_{2}\right) \right] \right\Vert _{X_{p}}\leq \left( M+1\right) \bar{f}%
\left( M+1\right) \left\Vert u_{1}-u_{2}\right\Vert _{X_{\infty }}\left\Vert
D_{\xi }^{2}u_{1}\right\Vert _{Y^{2,p}}^{2}+ 
\]%
\[
\bar{f}\left( M+1\right) \left\Vert D_{\xi }\left( u_{1}-u_{2}\right)
\right\Vert _{Y^{2,p}}\left\Vert D_{\xi }\left( u_{1}+u_{2}\right)
\right\Vert _{Y^{2,p}}+ 
\]%
\[
\bar{f}\left( M+1\right) \left\Vert u_{1}-u_{2}\right\Vert _{X_{\infty
}}\left\Vert D_{\xi }^{2}u_{1}\right\Vert _{X_{p}}+\bar{f}\left( M+1\right)
\left\Vert D_{\xi }\left( u_{1}-u_{2}\right) \right\Vert _{X_{p}}\leq 
\]%
\begin{equation}
C^{2}\bar{f}\left( M+1\right) \left\Vert u_{1}-u_{2}\right\Vert _{X_{\infty
}}\left\Vert u_{1}\right\Vert _{X_{\infty }}\left\Vert D_{\xi
}^{2}u_{1}\right\Vert _{X_{p}}+  \tag{3.13}
\end{equation}%
\[
C^{2}\bar{f}\left( M+1\right) \left\Vert u_{1}-u_{2}\right\Vert _{X_{\infty
}}^{\frac{1}{2}}\left\Vert D_{\xi }^{2}\left( u_{1}-u_{2}\right) \right\Vert
_{X_{p}}\left\Vert u_{1}+u_{2}\right\Vert _{X_{\infty }}^{\frac{1}{2}%
}\left\Vert D_{\xi }^{2}\left( u_{1}+u_{2}\right) \right\Vert _{X_{p}} 
\]%
\[
+\left( M+1\right) \bar{f}\left( M+1\right) \left\Vert
u_{1}-u_{2}\right\Vert _{X_{\infty }}+\bar{f}\left( M+1\right) \left\Vert
D_{\xi }^{2}\left( u_{1}-u_{2}\right) \right\Vert _{X_{p}}\leq 
\]%
\[
3C^{2}\left( M+1\right) ^{2}\bar{f}\left( M+1\right) \left\Vert
u_{1}-u_{2}\right\Vert _{X_{\infty }}+2C^{2}\left( M+1\right) \bar{f}\left(
M+1\right) \left\Vert D_{\xi }^{2}\left( u_{1}-u_{2}\right) \right\Vert
_{X_{p}}, 
\]%
where $C$ is the constant in Lemma $3.1$. From $\left( 3.10\right) -\left(
3.11\right) $, using Minkowski's inequality for integrals, Fourier
multiplier theorems for operator-valued functions in $X_{p}$ spaces and
Young's inequality, we obtain%
\[
\left\Vert G\left( u_{1}\right) -G\left( u_{2}\right) \right\Vert _{Y\left(
T\right) }\leq \dint\limits_{0}^{t}\left\Vert u_{1}-u_{2}\right\Vert
_{X_{\infty }}d\tau +\dint\limits_{0}^{t}\left\Vert u_{1}-u_{2}\right\Vert
_{Y^{2,p}}d\tau + 
\]%
\[
\dint\limits_{0}^{t}\left\Vert f\left( u_{1}\right) -f\left( u_{2}\right)
\right\Vert _{X_{\infty }}d\tau +\dint\limits_{0}^{t}\left\Vert f\left(
u_{1}\right) -f\left( u_{2}\right) \right\Vert _{Y^{2,p}}d\tau \leq 
\]%
\[
T\left[ 1+C_{1}\left( M+1\right) ^{2}\bar{f}\left( M+1\right) \right]
\left\Vert u_{1}-u_{2}\right\Vert _{Y\left( T\right) }, 
\]%
where $C_{1}$ is a constant. If $T$ satisfies $\left( 3.9\right) $ and the
following inequality 
\begin{equation}
T\leq \frac{1}{2}\left[ 1+C_{1}\left( M+1\right) ^{2}\bar{f}\left(
M+1\right) \right] ^{-1},  \tag{3.14}
\end{equation}%
then 
\[
\left\Vert Gu_{1}-Gu_{2}\right\Vert _{Y\left( T\right) }\leq \frac{1}{2}%
\left\Vert u_{1}-u_{2}\right\Vert _{Y\left( T\right) }. 
\]

That is, $G$ is a constructive map. By contraction mapping principle we know
that $G(u)$ has a fixed point $u(x,t)\in $ $Q\left( M;T\right) $ that is a
solution of the problem $(1.1)-(1.2)$. From $\left( 2.5\right) $ we get that 
$u$ is a solution of the following integral equation 
\[
u\left( t,x\right) =S_{1}\left( t,A\right) \varphi \left( x\right)
+S_{2}\left( t,A\right) \psi \left( x\right) + 
\]%
\[
+\dint\limits_{0}^{t}F^{-1}S\left( t-\tau ,\xi ,A\right) \left[ 1+L\left(
\xi \right) \right] ^{-1}\hat{f}\left( u\right) \left( \xi ,\tau \right)
d\tau ,\text{ }t\in \left( 0,T\right) . 
\]

Let us show that this solution is a unique in $Y\left( T\right) $. Let $%
u_{1} $, $u_{2}\in Y\left( T\right) $ are two solution of the problem $%
(1.1)-(1.2)$. Then%
\begin{equation}
u_{1}-u_{2}=\dint\limits_{0}^{t}F^{-1}S\left( t-\tau ,\xi ,A\right) \left[
1+L\left( \xi \right) \right] ^{-1}\left[ \hat{f}\left( u_{1}\right) \left(
\xi ,\tau \right) -\hat{f}\left( u_{2}\right) \left( \xi ,\tau \right) %
\right] d\tau .  \tag{3.15}
\end{equation}%
By definition of the space $Y\left( T\right) $, we can assume that%
\[
\left\Vert u_{1}\right\Vert _{X_{\infty }}\leq C_{1}\left( T\right) ,\text{ }%
\left\Vert u_{1}\right\Vert _{X_{\infty }}\leq C_{1}\left( T\right) . 
\]%
Hence, by Lemmas 2.3, Minkowski's inequality for integrals and Theorem 2.1
we obtain from $\left( 3.15\right) $

\begin{equation}
\left\Vert u_{1}-u_{2}\right\Vert _{Y^{2,p}}\leq C_{2}\left( T\right) \text{ 
}\dint\limits_{0}^{t}\left\Vert u_{1}-u_{2}\right\Vert _{Y^{2,p}}d\tau . 
\tag{3.16}
\end{equation}%
From $(3.16)$ and Gronwall's inequality, we have $\left\Vert
u_{1}-u_{2}\right\Vert _{Y^{2,p}}=0$, i.e. problem $(1.1)-(1.2)$ has a
unique solution which belongs to $Y\left( T\right) .$ That is, we obtain the
first part of the assertion. Now, let $\left[ 0\right. ,\left. T_{0}\right) $
be the maximal time interval of existence for $u\in Y\left( T_{0}\right) $.
It remains only to show that if $(3.3)$ is satisfied, then $T_{0}=\infty $.
Assume contrary that, $\left( 3.3\right) $ holds and $T_{0}<\infty .$ For $%
T\in \left[ 0\right. ,\left. T_{0}\right) $ we consider the following
integral equation

\bigskip 
\begin{equation}
\upsilon \left( x,t\right) =S_{1}\left( t,A\right) u\left( x,T\right)
+S_{2}\left( t,A\right) u_{t}\left( x,T\right) +  \tag{3.17}
\end{equation}

\[
\dint\limits_{0}^{t}F^{-1}S\left( t-\tau ,\xi ,A\right) \left[ 1+L\left( \xi
\right) \right] ^{-1}\hat{f}\left( \upsilon \right) \left( \xi ,\tau \right)
d\tau ,\text{ }t\in \left( 0,T\right) . 
\]%
By virtue of $(3.3)$, for $T^{\prime }>T$ we have 
\[
\sup_{t\in \left[ 0\right. ,\left. T\right) }\left( \left\Vert u\right\Vert
_{Y^{2,p}}+\left\Vert u\right\Vert _{X_{\infty }}+\left\Vert
u_{t}\right\Vert _{Y^{2,p}}+\left\Vert u_{t}\right\Vert _{X_{\infty
}}\right) <\infty . 
\]

By reasoning as a first part of theorem and by contraction mapping
principle, there is a $T^{\ast }\in \left( 0,T_{0}\right) $ such that for
each $T\in \left[ 0\right. ,\left. T_{0}\right) ,$ the equation $\left(
3.17\right) $ has a unique solution $\upsilon \in Y\left( T^{\ast }\right) .$
The estimates $\left( 3.9\right) $ and $\left( 3.14\right) $ imply that $%
T^{\ast }$ can be selected independently of $T\in \left[ 0\right. ,\left.
T_{0}\right) .$ Set $T=T_{0}-\frac{T^{\ast }}{2}$ and define 
\begin{equation}
\tilde{u}\left( x,t\right) =\left\{ 
\begin{array}{c}
u\left( x,t\right) ,\text{ }t\in \left[ 0,T\right] \\ 
\upsilon \left( x,t-T\right) \text{, }t\in \left[ T,T_{0}+\frac{T^{\ast }}{2}%
\right]%
\end{array}%
\right. .  \tag{3.18}
\end{equation}

By construction $\tilde{u}\left( x,t\right) $ is a solution of the problem $%
(1.1)-(1.2)$ on $\left[ T,T_{0}+\frac{T^{\ast }}{2}\right] $ and in view of
local uniqueness, $\tilde{u}\left( x,t\right) $ extends $u.$ This is against
to the maximality of $\left[ 0\right. ,\left. T_{0}\right) $, i.e we obtain $%
T_{0}=\infty .$\qquad

\begin{center}
\bigskip \textbf{4. \ The Cauchy problem for the system of Boussinesq
equation }
\end{center}

\bigskip Consider the Cauchy problem for the following nonlinear system%
\begin{equation}
\left( u_{m}\right) _{tt}-\left( Lu_{m}\right)
_{tt}+\sum\limits_{j=1}^{N}a_{mj}u_{j}\left( x,t\right) =f_{m}\left(
u\right) ,\text{ }x\in R^{n},\text{ }t\in \left( 0,T\right) ,  \tag{4.1}
\end{equation}%
\begin{equation}
u_{m}\left( x,0\right) =\varphi _{m}\left( x\right) ,\text{ }\frac{\partial 
}{\partial t}u_{m}\left( x,0\right) =\psi _{m}\left( x\right) ,\text{ }%
m=1,2,...,N,\text{ }N\in \mathbb{N},  \tag{4.2}
\end{equation}%
where $u=\left( u_{1},u_{2},...,u_{N}\right) ,$ $a_{mj}$ are complex
numbers, $\varphi _{m}\left( x\right) $, $\psi _{m}\left( x\right) $ are
data functions and 
\[
Lu=\dsum\limits_{i,j=1}^{2}a_{ij}\frac{\partial ^{2}u}{\partial
x_{i}\partial x_{j}},\text{ }a_{ij}\in \mathbb{C}. 
\]%
Let%
\[
l_{q}=\text{ }l_{q}\left( N\right) =\left\{ \text{ }u=\left\{ u_{j}\right\} ,%
\text{ }j=1,2,...N,\left\Vert u\right\Vert _{l_{q}\left( N\right) }=\left(
\sum\limits_{j=1}^{N}\left\vert u_{j}\right\vert ^{q}\right) ^{\frac{1}{q}%
}<\infty \right\} , 
\]%
(see $\left[ \text{44, \S\ 1.18}\right] .$ Let $A$ be the operator in $%
l_{q}\left( N\right) $ defined by%
\[
\text{ }A=\left[ a_{mj}\right] \text{, }a_{mj}=g_{m}2^{sj},\text{ }%
m,j=1,2,...,N,\text{ }D\left( A\right) =\text{ }l_{q}^{s}\left( N\right) = 
\]

\[
\left\{ \text{ }u=\left\{ u_{j}\right\} ,\text{ }j=1,2,...N,\left\Vert
u\right\Vert _{l_{q}^{s}\left( N\right) }=\left(
\sum\limits_{j=1}^{N}2^{sj}u_{j}^{q}\right) ^{\frac{1}{q}}<\infty \right\} . 
\]

Let 
\[
X_{pq}=L^{p}\left( R^{n};l_{q}\right) ,Y^{s,p,q}=L^{s,p}\left(
R^{n};l_{q}\right) ,Y_{1}^{s,p,q}=L^{s,p}\left( R^{n};l_{q}\right) \cap
L^{1}\left( R^{n};l_{q}\right) , 
\]

\[
Y_{\infty }^{s,p,q}=L^{s,p}\left( R^{n};l_{q}\right) \cap L^{\infty }\left(
R^{n};l_{q}\right) 
\]

and%
\[
\text{ }Y^{2,p,q}=W^{2,p}\left( R^{n};l_{q}^{s},l_{q}\right) ,\text{ }%
E_{0q}=B_{p,p}^{2\left( 1-\frac{1}{2p}\right) }\left( R^{n};\left(
l_{q}^{s},l_{q}\right) _{\frac{1}{2p},p},l_{q}\right) . 
\]

From Theorem 3.1 we obtain the following result

\textbf{Theorem 4.1. }Let Conditon 2.0 hold. Assume $\varphi _{m},$ $\psi
_{m}$ $\in Y_{\infty }^{2,p,q}$ and $1<p<\infty $ for$\frac{n}{p}<2$.
Suppose the function $u\rightarrow $ $f\left( u\right) $: $R^{n}\times \left[
0,T\right] \times E_{0q}\rightarrow l_{q}$ is a measurable function in $%
\left( x,t\right) \in R^{n}\times \left[ 0,T\right] $ for $u\in E_{0q};$ $%
f\left( x,t.,.\right) $ and this function is continuous in $u\in E_{0q}$ for 
$x,t\in R^{n}\times \left[ 0,T\right] ;$ moreover $f\left( u\right) \in
C^{\left( 3\right) }\left( E_{0q};l_{q}\right) $. Then problem $\left(
4.1\right) -\left( 4.2\right) $ has a unique local strange solution $u\in
C^{\left( 2\right) }\left( \left[ 0\right. ,\left. T_{0}\right) ;Y_{\infty
}^{2,p,q}\right) $, where $T_{0}$ is a maximal time interval that is
appropriately small relative to $M$. Moreover, if

\begin{equation}
\sup_{t\in \left[ 0\right. ,\left. T_{0}\right) }\left( \left\Vert
u\right\Vert _{Y^{2,p,q}}+\left\Vert u\right\Vert _{X_{\infty
,q}}+\left\Vert u_{t}\right\Vert _{Y^{2,p,q}}+\left\Vert u_{t}\right\Vert
_{X_{\infty ,q}}\right) <\infty  \tag{4.3}
\end{equation}%
then $T_{0}=\infty .$

\ \textbf{Proof. }By virtue of $\left[ 43\right] ,$ the $l_{q}\left(
N\right) $ is a UMD space. It is easy to see that the operator $A$ is $R$%
-positive in $l_{q}\left( N\right) .$ Moreover, by interpolation theory of
Banach spaces $\left[ \text{44, \S\ 1.3}\right] $, we have 
\[
E_{0q}=\left( W^{2,p}\left( R^{n};l_{q}^{s},l_{q}\right) ,L^{p}\left(
R^{n};l_{q}\right) \right) _{\frac{1}{2p},q}=B_{p,q}^{2\left( 1-\frac{1}{2p}%
\right) }\left( R^{n};l_{q}^{s\left( 1-\frac{1}{2p}\right) },l_{q}\right) . 
\]%
By using the properties of spaces $Y^{s,p,q},$ $Y_{\infty }^{s,p,q},$ $%
E_{0q} $ we get that all conditions of Theorem 3.1 are hold, i,e., we obtain
the conclusion.\ 

\begin{center}
\textbf{Acknowledgements}
\end{center}

The author would like to express a gratitude to Bulent Eryigit for his
useful advices in English in preparing of this paper.

\textbf{References}

\begin{quote}
\ \ \ \ \ \ \ \ \ \ \ \ \ \ \ \ \ \ \ \ \ \ \ \ 
\end{quote}

\begin{enumerate}
\item V.G. Makhankov, Dynamics of classical solutions (in non-integrable
systems), Phys. Lett. C 35, (1978) 1--128.

\item G.B. Whitham, Linear and Nonlinear Waves, Wiley--Interscience, New
York, 1975.

\item N.J. Zabusky, Nonlinear Partial Differential Equations, Academic
Press, New York, 1967.

\item C. G. Gal and A. Miranville, Uniform global attractors for
non-isothermal viscous and non-viscous Cahn--Hilliard equations with dynamic
boundary conditions, Nonlinear Analysis: Real World Applications 10 (2009)
1738--1766.

\item T. Kato, T. Nishida, A mathematical justification for Korteweg--de
Vries equation and Boussinesq equation of water surface waves, Osaka J.
Math. 23 (1986) 389--413.

\item A. Clarkson, R.J. LeVeque, R. Saxton, Solitary-wave interactions in
elastic rods, Stud. Appl. Math. 75 (1986) 95--122.

\item P. Rosenau, Dynamics of nonlinear mass-spring chains near continuum
limit, Phys. Lett. 118A (1986) 222--227.

\item S. Wang, G. Chen, Small amplitude solutions of the generalized IMBq
equation, J. Math. Anal. Appl. 274 (2002) 846--866.

\item S. Wang, G. Chen,The Cauchy Problem for the Generalized IMBq equation
in $W^{s,p}\left( R^{n}\right) $, J. Math. Anal. Appl. 266, 38--54 (2002).

\item Denk R., Hieber M., Pr\"{u}ss J., $R$-boundedness, Fourier multipliers
and problems of elliptic and parabolic type, Mem. Amer. Math. Soc. 166
(2003), n.788.

\item H. O. Fattorini, Second order linear differential equations in Banach
spaces, in North Holland Mathematics Studies,\ V. 108, North-Holland,
Amsterdam, 1985.

\item J. A. Goldstein, Semigroup of linear operators and applications,
Oxford, 1985.

\item M. Sova, \textquotedblleft Linear differential equations in Banach
spaces,\textquotedblright\ Rozpr. CSAV MPV, 85 (6), 1--86 (1975)

\item G. Da Prato and E. Giusti, A characterization on generators of
abstract cosine functions,\ Boll. Del. Unione Mat., (22)1967, 367--362.

\item A. Arosio, Linear second order differential equations in Hilbert
space. The Cauchy problem and asymptotic behaviour for large time, Arch.
Rational Mech. Anal., v 86, (2), 1984, 147-180.

\item N. Tanaka and I. Miyadera, Exponentially bounded C-semigroups and
integrated semi-groups, Tokyo J. Math. (12) 1989, 99-115.

\item E. B. Davies and M. M. Pang, The Cauchy problem and a generalization
of the Hille-Yosida theorem, Proc. London Math. Soc. (55) 1987, 181-208.

\item S. Piskarev and S.-Y. Shaw, Multiplicative perturbations of semigroups
and applications to step responses and cumulative outputs, J. Funct. Anal.
128 (1995), 315-340.

\item S. Piskarev and S.-Y. Shaw, Perturbation and comparison of cosine
operator functions, Semigroup Forum 51 (1995), 225-246.

\item S. Piskarev and S.-Y. Shaw, On certain operator families related to
cosine operator functions,Taiwan. J. Math., (1)1997, 527--546.

\item H. Amann, Linear and quasi-linear equations,1, Birkhauser, Basel
1995.\ \ \ \ \ \ \ \ \ \ \ \ \ \ \ \ \ \ \ \ \ \ \ \ \ \ \ \ \ \ \ \ \ \ \ \
\ \ \ \ \ \ \ \ \ \ \ \ \ \ \ \ \ \ \ \ \ \ \ \ \ \ \ 

\item R. Agarwal R, D. O' Regan, V. B. Shakhmurov, Separable anisotropic
differential operators in weighted abstract spaces and applications, J.
Math. Anal. Appl. 2008, 338, 970-983.

\item A. Ashyralyev, Claudio Cuevas and S. Piskarev, On well-posedness of
difference schemes for abstract elliptic problems in spaces, Numer. Func.
Anal.Opt., v. 29, No. 1-2, Jan. 2008, 43-65

\item A. Favini, V. Shakhmurov, Y. Yakubov, Regular boundary value problems
for complete second order elliptic differential-operator equations in UMD
Banach spaces, Semigroup form, v. 79 (1), 2009.

\item V. I. Gorbachuk and M. L. Gorbachuk M, Boundary value problems \ for
differential-operator equations, Naukova Dumka, Kiev, 1984.

\item S. G. Krein, Linear differential equations in Banach space,
Providence, 1971.

\item A. Lunardi, Analytic semigroups and optimal regularity in parabolic
problems, Birkhauser, 2003.

\item J. L Lions and J. Peetre, Sur une classe d'espaces d'interpolation,
Inst. Hautes Etudes Sci. Publ. Math., 19(1964), 5-68.

\item P. E. Sobolevskii, Inequalities coerciveness for abstract parabolic
equations, Dokl. Akad. Nauk. SSSR, (1964), 57(1), 27-40.

\item A.Ya. Shklyar, Complete second order linear differential equations in
Hilbert spaces, Birkhauser Verlag, Basel, 1997.

\item V. B. Shakhmurov, Embedding and maximal regular differential operators
in Banach-valued weighted spaces\textbf{, }Acta mathematica Sinica, 22(5)
(2006), 1493-1508.

\item V. B. Shakhmurov, Embedding and separable differential operators in
Sobolev-Lions type spaces, Mathematical Notes, 2008, v. 84, no 6, 906-926.

\item V. B. Shakhmurov, Nonlinear abstract boundary value problems in
vector-valued function spaces and applications, Nonlinear Analysis Series A:
Theory, Method \& Applications, v. 67(3) 2006, 745-762.

\item V. B. Shakhmurov, Embedding theorems and\ maximal regular differential
operator equations in Banach-valued function spaces, Journal of Inequalities
and Applications, 2( 4), (2005), 329-345.

\item V. B. Shakhmurov, The Cauchy problem for generalized abstract
Boussinesq equations, Dynamic systems and applications, 25, (2016),109-122.

\item V. B. Shakhmurov, Nonlocal problems for Boussinesq equations,
Nonlinear Analysis TMA, 142 (2016) 134-151.

\item V. B. Shakhmurov, Separable anisotropic differential operators and
applications, J. Math. Anal. Appl. 2006, 327(2), 1182-1201.

\item P. Guidotti, \ Optimal regularity for a class of singular abstract
parabolic equations, J. Differ. Equations, v. 232, 2007, 468--486.

\item R. Shahmurov, On strong solutions of a Robin problem modeling heat
conduction in materials with corroded boundary, Nonlinear Analysis,
Real-wold applications, v.13, (1), 2011, 441-451.

\item R. Shahmurov, Solution of the Dirichlet and Neumann problems for a
modified Helmholtz equation in Besov spaces on an annuals, Journal of
differential equations, v. 249(3), 2010, 526-550.

\item L. Weis, Operator-valued Fourier multiplier theorems and maximal $%
L_{p} $ regularity, Math. Ann. 319, (2001), 735-758.

\item S. Yakubov and Ya. Yakubov, \textquotedblright Differential-operator
Equations. Ordinary and Partial \ Differential Equations \textquotedblright
, Chapman and Hall /CRC, Boca Raton, 2000.

\item D. L. Burkholder: A geometric condition that implies the existence of
certain singular integrals of Banach-space-valued functions, in: Conference
on harmonic analysis in honor of Antoni Zygmund, V. I, II, Chicago,
Illinois, 1981, 270-286.

\item H. Triebel, Interpolation theory, Function spaces, Differential
operators, North-Holland, Amsterdam, 1978.

\item H. Triebel, Fractals and spectra, Birkhauser Verlag, Related to
Fourier analysis and function spaces, Basel, 1997.

\item R. Haller, H. Heck, A. Noll, Mikhlin's theorem for operator-valued
Fourier multipliers in $n$ variables, Math. Nachr. 244 (2002), 110-130.

\item L. Nirenberg, On elliptic partial differential equations, Ann. Scuola
Norm. Sup. Pisa 13 (1959), 115--162.

\item S. Klainerman, Global existence for nonlinear wave equations, Comm.
Pure Appl. Math. 33 (1980), 43--101.

\item J. L. Lions and J. Peetre, Sur une classe d'espaces d'interpolation,
Inst. Hautes Etudes Sci. Publ. Math., 19(1964), 5-68.

\item H. O. Fattorini, \textquotedblleft Ordinary differential equations in
linear topological spaces, I,\textquotedblright\ J. Differential Equations,
5, No. 1, 72--105 (1969).

\item C. Travis and G. F. Webb, Second order differential equations in
Banach spaces, Nonlinear Equations in Abstract Spaces, (ed. by V.
Lakshmikantham), Academic Press, 1978, 331-361.

\item G. Da Prato and E. Giusti, Una caratterizzazione del generator di
funzioni coseno astratte, Boll. Un. Mat. Ital., 22 (1967), 357-362.

\item Y. Xiao and Z. Xin, On the vanishing viscosity limit for the 3D
Navier-Stokes equations with a slip boundary condition. Comm. Pure Appl.
Math. 60, 7 (2007), 1027--1055.
\end{enumerate}

\end{document}